\documentclass[conference]{IEEEtran}
\IEEEoverridecommandlockouts
\usepackage{amsmath,amssymb,amsfonts}
\usepackage{algorithmic}
\usepackage{graphicx}
\usepackage{textcomp}
\usepackage{xcolor}
\usepackage{comment}
\def\BibTeX{{\rm B\kern-.05em{\sc i\kern-.025em b}\kern-.08emT\kern-.1667em\lower.7ex\hbox{E}\kern-.125emX}}
\usepackage[
backend=biber,
style=numeric,
sorting=none
]{biblatex}

\begin{document}
\title{A novel DC Microgrid stabilization method based on the logarithmic norm\\}

\thanks{This work was partially supported by the Prime Minister's Research Fellowship of the Government of India.
}

\author{\IEEEauthorblockN{Aditya Arpitha Prasad}
\IEEEauthorblockA{\textit{Computational \& Data Sciences} \\
\textit{Indian Institute of Science}\\
\textit{Bengaluru 560012, India} \\
\textit{adityaprasad@iisc.ac.in}} 
}

\maketitle

\begin{abstract}
Microgrids could be the answer to integrating distributed energy resources into our powergrid. It promises improved resilience, reliability, efficiency and decarbonising of our electric grid. This paper models a low voltage direct current microgrid using a DAE formulation on a complex network. We then offer a novel method to compute a stability metric using the logarithmic norm which can be used to isolate disturbed nodes.
\end{abstract}


\section{Introduction}
There are well documented reasons motivating a shift towards a decentralized grid incorporating more renewable energy sources - security \cite{Hirsch_Parag_Guerrero_2018}, economic \cite{Ahmad_Alam_2018}, ecological \cite{change2014mitigation}, social \cite{Jhunjhunwala_Kaur_2018}. A decentralized grid is more resistant to natural disasters, cyber attacks and safeguard our ability to secure uninterrupted quality power to critical infrastructures like hospitals, power plants, network towers. Due to the ability of individual microgrids to island itself, the grid will be resistant to cascading blackouts, improving resilience and reliability.

On the ecological side, the energy supply sector is the largest contributor to global greenhouse gas emissions and decarbonizing electricity generation is a cost-effective way of rapidly achieving our carbon targets. As energy generation comes closer to the point of consumption, we can bring down line losses. Consumers will turn into prosumers with bidirectional power flow on the grid. We can also adopt efficient practices of buffering energy for peak demand instead of the costly inefficient ramp up that now occurs. For a country like India, where many communities live in remote and inaccessible locations microgrids can help  bridging last mile connectivity \cite{Heynen_Lant_Smart_Sridharan_Greig_2019}.

India has set a target of 175 GW installed capacity by 2022 of which 100 GW is solar \cite{Palchak_Cochran_Deshmukh_Ehlen_Soonee_Narasimhan_Joshi_McBennett_Milligan_Sreedharan_et_al_2017}. These are ambitious goals and microgrids are a natural way to incorporate these  distributed energy resources (DER) in a decentralized way. But there are challenges to overcome when integrating DERs into the main grid. Due to cloud cover and wind speed variability the grid has to handle both supply and demand side uncertainties \cite{Steen_Goop_Göransson_Nursbo_Brolin_2014}. As renewable penetration increases we see a reduction in inertia. There is also a need for research into new distributed control techniques to ensure robust control of the microgrid.

There are two main types of microgrids - AC and DC. Due to the prevalence of DC loads, DC sources and reduced losses when compared to AC-DC power electronics, DC microgrid is what we used in the model. DC microgrids are characterized by components such as controlled and uncontrolled loads, distributed generation (DG) units and storage devices operate together in a coordinated manner with controlled power electronic devices under a DC power system principle \cite{Justo_Mwasilu_Lee_Jung_2013}. Even in case of sources like wind energy, the AC power that is produced is variable and needs to be rectified to DC then inverted to be used in AC grids. EPRI found that if we switched to small DC grids we could reduce losses by 10-30\% \cite{Patterson_2012} 

At present research on stability of microgrids in islanded mode and grid connected mode are treated separately, here we are dealing with islanded mode since the utility grid has sufficient inertia to handle any disturbances in very short time scales \cite{Shuai_Sun_Shen_Tian_Tu_Li_Yin_2016}. Further there are two types of small signal stability or transient stability. Most of this analysis is done on extremely small number of DGs under consideration, here we will use the logarithmic norm to consider the local stability at the equilibrium point.

Recent work usually deals with a linearized the system and exploits Lypunov function to prove local stability \cite{Belk_Inam_Perreault_Turitsyn_2016}. Due to practical limitations in communication links, instead of a centralized control scheme like master-slave control, decentralized control schemes were explored. In most studies droop control is what is used as feedback control. \cite{Dag_Mirafzal_2016}

The traditional method to assess stability is eigenvalue analysis\cite{Renjit_Illindala_Yedavalli_2014}\cite{Su_Liu_Sun_Han_Hou_2018}\cite{Tamimi_Cañizares_Bhattacharya_2013}. This is without accounting for the fact that numerical errors and limitations in measurement can mean pseudospectrum of the matrix is often more relevant to the study \cite{Liu_Zhang_Rizzoni_2017}. The logarithmic norm gives us an estimate for psudeospectra \cite{Trefethen_Embree_2005} and stiffness \cite{Higham_Trefethen_1993}.

\section{Microgrid Model}

In this section we present the models used for the power lines, distributed sources, constant power load and complex network data structure used to represent a DC microgrid.

\subsection{Topology}

We can describe the network as a weighted, undirected graph ($\mathcal{N},\mathcal{E}$), where we have $n = |\mathcal{N}|$ nodes and $m = |\mathcal{E}|$ edges. The topology of this graph is defined by an incidence matrix $L_{m \times n}$, whenever a new node i is connected to the DC microgrid at bus/node j a new row is appended to this matrix where $L(m,i) = 1$ and $L(m,j) = -1$. This lets us use $L \in \mathbb R^{m \times n}$ to set up the KCL and KVL constraints easily by accounting for current direction as shown below in equations (1) and (2). Here let $\vec v_g(t) = [v_{g1}, v_{g2}, v_{g3}, \cdots, v_{gn}]$ be the bus voltages. Similarly $\vec v_p(t) = [v_{p1}, v_{p2}, v_{p3}, \cdots, v_{pm}]$ are the power line voltages found as the difference between node voltages.

\begin{equation}
    \vec v_p(t)_{m \times 1} = L_{m \times n}\vec v_g(t)_{n \times 1}
\end{equation}

For the current law we can define $\vec i_g(t) = [i_{g1}, i_{g2}, i_{g3}, \cdots, i_{gn}]$ be the current leaving each node. Similarly $\vec i_p(t) = [i_{p1}, i_{p2}, i_{p3}, \cdots, i_{pm}]$ are the power line currents found flowing 

\begin{equation}
    \vec i_g(t)_{n \times 1} = L^T_{n \times m}\vec i_p(t)_{n \times 1}
\end{equation}

\subsection{Power line model}

Since this a DC microgrid with short line lengths, we can go for a nominal pi model of a transmission line, where we have a lumped resistance and inductance in series. This is a common model in the literature \cite{Anand_Fernandes_2013}. This gives us the differential equation,

\begin{equation}
   \vec v_p(t) = I_{n \times n}\frac{d \vec i_p(t)}{dt} + R_{n \times n}\vec i_p(t)
\end{equation}

Here the matrix $I$ is a diagonal matrix of the line inductances and $R$ is the diagonal matrix of the line resistances. 

\begin{figure}[htbp]
\centerline{\includegraphics[width=0.5\textwidth]{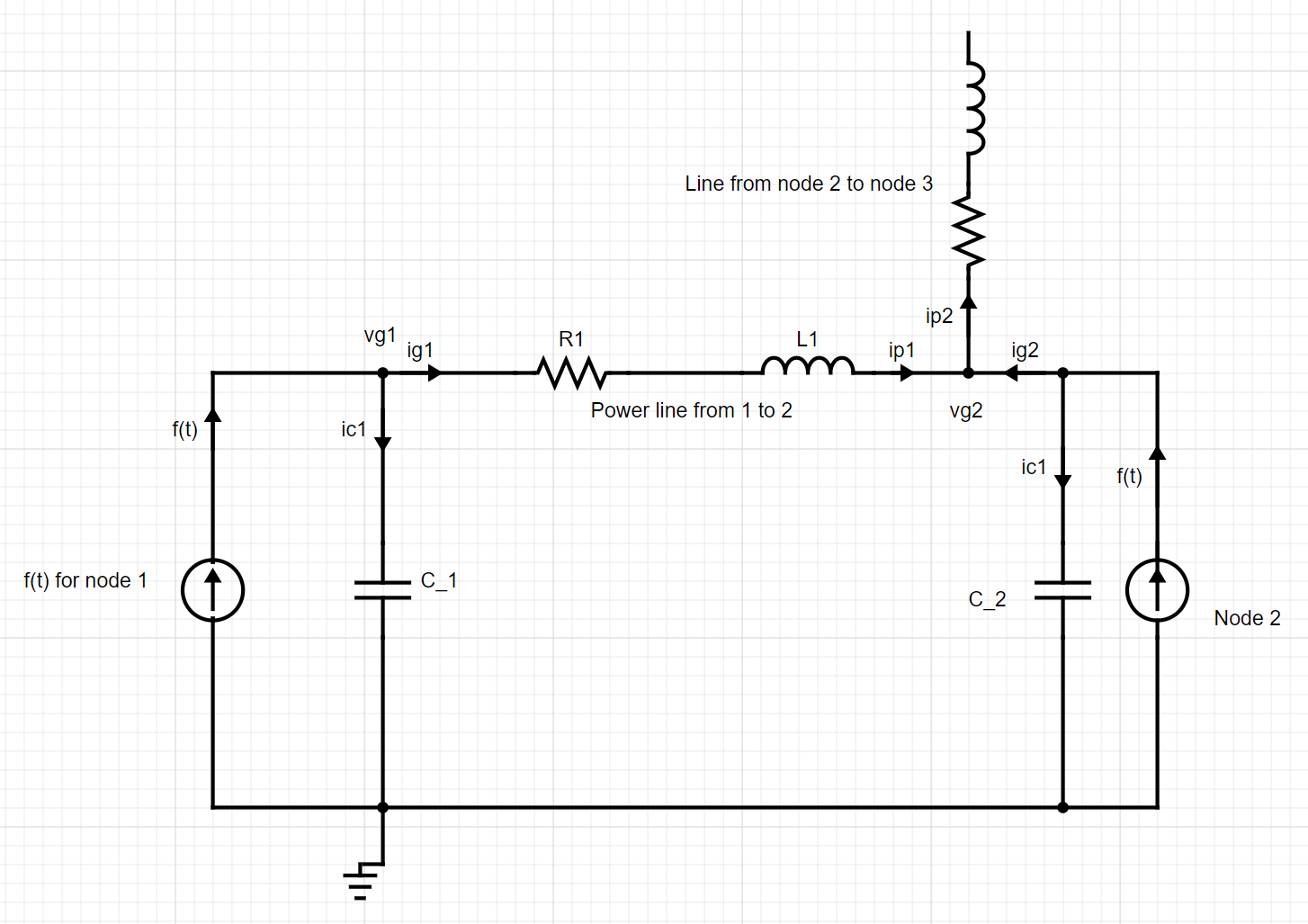}}
\caption{Part of a DC micro grid with voltages and currents labelled}
\label{fig}
\end{figure}

\subsection{PV Solar Home}

Solar photovoltaic (PV) systems can be set up in each home letting it either pull power or supply power to the microgrid. The DC microgrid will be able to grow via preferential attachment organically. When a node is acting as a net producer we can assume a controlled current source and capacitance in parallel to model the DC link capacitance as shown in Fig 1. This is because modern power electronics can supply that constant current after a DC-DC step up from an unregulated dc voltage. \cite{Ansari_Chandel_Tariq_2021}

We can define the current entering each node's capacitance as $\vec i_c(t) = [i_{c1}, i_{c2}, i_{c3}, \cdots, i_{cn}]$. $C$ is a diagonal matrix with all the link capacitance values. 

\begin{equation}
    \vec i_g(t) = \vec f(t) - \vec i_c(t)
\end{equation}

\begin{equation}
    \vec i_c(t) = C_{n \times n} \frac{d \vec v_g(t)}{dt}
\end{equation}

\subsection{DAE formulation}

From these equations we can formulate the inherent ODE of the DC microgrid, from plugging equations (1) into (3),

\begin{equation}
    \frac{d \vec i_p(t)}{dt}  = I^{-1} ( L \vec v_g(t) - R \vec i_p(t) )
\end{equation}

Now by plugging (4) and (2) in (5)

\begin{equation}
    \frac{d \vec v_g(t)}{dt}  = C^{-1} ( -L^T \vec i_p(t) + f(t) )
\end{equation}

Now by defining the state vector as $x = [ i_p(t), v_g(t) ] $,

\begin{equation}
    \dot x =  \begin{bmatrix}
    -I^{-1}R && I^{-1}L  \\
    -C^{-1}L^T && O 
\end{bmatrix} x + \begin{bmatrix}
    0 \\
    C^{-1} 
\end{bmatrix} f(t)
\end{equation}

\section{Stability Analysis}

\subsection{Proportional Droop Control}

For my analysis I will be using proportional droop control where the control law is based on the local measured $v_g$ value as compared to a set reference value - $V_s$. The matrix $K$ For these producer nodes the current source will be controlled to be. We can model the consuming nodes as a constant power load - $P_s$ since we can assume the power electronics tightly regulates the output \cite{Zhang_Yan_2011}.

Here the model of $f$ will change depending on which nodes are producing excess power and supplying to the grid and which are drawing current as a load. I use a permutation matrix $P$ to order the $p$ power producing nodes in the front and $n-p$ nodes will be in the end of the vector.

\begin{equation}
    f_p(t) = \hat K_{p \times p} (\hat v_g(t) - V_s(t) )
\end{equation}

\begin{equation}
    f_c(t) = \frac{P_s}{V_s}
\end{equation}

Where the hat is used to represent the transformed variables via permutation. 

$$\hat v_g = P \vec v_g$$
$$\hat i_g = P \vec i_g$$
$$\hat i_c = P \vec i_c$$
$$\hat L = P L$$
$$\hat R = P R$$
$$\hat I = P I$$
$$\hat C = P C$$

We can write the complete $f$ as,

\begin{equation}
    f(t) = \begin{bmatrix}
    K_{p \times p} && 0 \\
    0 && 0 
\end{bmatrix} \hat v_g - \begin{bmatrix}
    K_{p \times p} && 0 \\
    0 && 0 
\end{bmatrix}  \begin{bmatrix}
    V_s\\
    -\frac{P_s}{V_s} 
\end{bmatrix}
\end{equation}

So the sake of conciseness we can define the constant matrix $K'$ and vector $c$ using equation (8),

\begin{equation}
    f(t) = K'\hat v_g - K'c
\end{equation}
 
So using equations (9) and (10) in (7),

\begin{equation}
    \frac{d \hat v_g(t)}{dt}  = \hat C^{-1} \left ( -\hat L^T i_p(t) +  K'\hat v_g - K'c \right )
\end{equation}

Now by redefining the state vector as $z = [ \vec i_p(t), \hat v_g(t) ] $, we can use this to formulate the closed loop system,

\begin{equation}
    \dot z =  \begin{bmatrix}
    -\hat I^{-1}\hat R && \hat I^{-1}L  \\
    -\hat C^{-1}\hat L^T && \hat C^{-1}K' 
\end{bmatrix} z + \begin{bmatrix}
    0 \\
    -C^{-1}K'c 
\end{bmatrix}
\end{equation}

For sake of consciseness again we define $B$ and $k$ using equation (14)

\begin{equation}
    \dot z = Bz + k
\end{equation}

This requires computation of inverse of a matrix where the numerical stability is not assured. But given a permutation matrix we can find the matrix $B$.

\subsection{Logarithmic norm}

We propose it can be useful to use the logarithmic norm introduced by Dahlquist \cite{Schmidt_1961} to measure the stability of the microgrid for any given permutation matrix $P$,

\begin{equation}
    \mu [B] = \lim_{h \to 0^+}\frac{||I+hB||-1}{h}
\end{equation}

Where the matrix norm used is the matrix norm induced by the vector 2 - norm. 

\begin{equation}
    ||B|| = \sup_{x\neq 0} \frac{||Bx||_2}{||x||_2}
\end{equation}

One important reason to prefer the logarithmic norm is because it captures the information regarding the pseudospectrum of $A$. Due to humidity, temperature, numerical errors, measurement errors we cannot be certain regarding the local stability of equilibrium of $B$ matrix in equation (15) from just the eigenvalues. The = $\epsilon$-pseudospectral abscissa and logarithmic norm, also called numerical abscissa are intimately connected. The logarithmic norm is maximum real part of the numerical range.\cite{Trefethen_Embree_2005} This is of a serious concern when it comes to nonnormal system matrices where even a small amount of noise can destabilize when eigenvalues are negative. \cite{Higham_Mao_2006}

\begin{equation}
    \mu [B] = \sup_{\epsilon > 0} [\alpha_\epsilon(B) - \epsilon]
\end{equation}

There are other useful properties like $\mu [B]$ is the maximal growth rate of $\log ||z||$ where $\dot z = Bz$. 
\begin{equation}
    D_t^+\log ||x|| \leq \mu [B] 
\end{equation}

Here $D_t^+$ is the upper right Dini derivative. 

So from equation (15) we can write the exponential stability condition $\mu [B] < 0$ \cite{Söderlind_2006}.

Here the computation of the logarithmic norm can be reduced to an eigenvalue problem, \cite{Mengi_Overton_2005}

\begin{equation}
    \mu [B] = \lambda_{max} \{\frac{B + B^*}{2}\}
\end{equation}

Simple eigenvalue computations alone cannot reveal the potential for transient growth. The spectral abscissa gives us a lower bound for the solutions, while the logarithmic norm gives us the upper bound. \cite{Perov_Kostrub_2017}

\begin{equation}
    e^{t\alpha} \leq ||e^{tB}|| \leq e^{t\mu[B]}
\end{equation}

Where $\alpha$ is the spectral abscissa, 
$$\alpha = \max_{1\leq i \leq n } \Re\{\lambda_i(B)\}$$

\subsection{Stabilization method}

Step 1: Monitor the deviation from voltage set point and if positive threshold $t$ is crossed, move the node to producer list and if negative threshold is crossed, move node to consumer list.

Step 2: Computer logarithmic norm for current B and store in a temporary variable.

Step 3: Recompute the matrix P with new producer and consumer list and recalculate the system matrix B.

Step 4: Calculate logarithmic norm for B and compare to old value.

Step 5: If value is more negative, turn switches to implement the move.

With this method the consumers and producers of the microgrid can be decided and stability can be iteratively improved.

\section*{Conclusion}

This paper has developed a complete model of a general scalable DC microgrid model and discussed a method to control switching between nodes. To the best of the author's knowledge this is the first time logarithmic norm has been applied to stability studies of DC microgrids. This paper proposed a novel technique to increase the stability of the microgrid. DC microgrids offers a lot of benefits and in future studies this model needs to be improved on to incorporate other DERs like wind farms, energy storage devices, fuel cells.

\printbibliography

\vspace{12pt}

\end{document}